%% file: LP.tex
\documentclass[12pt,reqno]{amsart}
\usepackage{amssymb,amscd}
\usepackage[all]{xy}
\usepackage[headings]{fullpage}
\usepackage{times}

\include{formatting}
\include{macros}

\newcommand{\Fsep}{F^{{sep}}}

\begin{document}
\title{Explicit points on the Legendre curve}
\author{Douglas Ulmer}
\address{School of Mathematics \\ Georgia Institute of Technology
  \\ Atlanta, GA 30332}
\email{ulmer@math.gatech.edu}
\thanks{This paper is based upon work partially supported by the National
Science Foundation under Grant No. DMS 0701053}
\subjclass[2010]{11G05, 11G40} 

\begin{abstract}
  We study the elliptic curve $E$ given by $y^2=x(x+1)(x+t)$ over the
  rational function field $k(t)$ and its extensions
  $K_d=k(\mu_d,t^{1/d})$.  When $k$ is finite of characteristic $p$
  and $d=p^f+1$, we write down explicit points on $E$ and show by
  elementary arguments that they generate a subgroup $V_d$ of rank
  $d-2$ and of finite index in $E(K_d)$.  Using more sophisticated
  methods, we then show that the Birch and Swinnerton-Dyer conjecture
  holds for $E$ over $K_d$, and we relate the index of $V_d$ in
  $E(K_d)$ to the order of the Tate-Shafarevich group $\sha(E/K_d)$.
  When $k$ has characteristic 0, we show that $E$ has rank 0 over
  $K_d$ for all $d$.
\end{abstract}

\maketitle

\section{Introduction}
Throughout the paper, $p>2$ is a prime number, $\Fp$ is the field of
$p$ elements, $\Fq$ is a finite extension of $\Fp$ of cardinality $q$,
and $\Fpbar$ is an algebraic closure of $\Fp$.  We write $\Fp(t)$,
$\Fq(t)$, and $\Fpbar(t)$ for the rational function fields over $\Fp$,
$\Fq$, and $\Fpbar$ respectively.

There are now several constructions of elliptic curves (and
higher-dimensional Jacobians) of large rank over $\Fp(t)$ or
$\Fpbar(t)$.  The first results in this direction are due to
Shafarevich and Tate \cite{TateShafarevich67}, and their arguments as
well as more recent results on large ranks are discussed in
\cite{Ulmer11}.  Most of these constructions rely on relatively
sophisticated mathematics, such as the theory of algebraic surfaces,
cohomology, and $L$-functions.

Our first aim in this note is to give an explicit, elementary
construction of points on a certain non-isotrivial elliptic curve over
$\Fp(t)$, the Legendre curve.  We then show that the rank of the group
these points generate can be arbitrarily large.  The arguments are
elementary and could have been given fifty years ago.\footnote{We
  note that in \cite{ConceicaoThesis}, Concei\c{c}\~ao gives an
  equally simple construction of polynomial points on certain {\it
    isotrivial\/} elliptic curves over $\Fp(t)$ and uses them to show
  that the rank is large.}  It is remarkable that such a classical and
well-studied object as the Legendre curve still has a few surprises.

Further pleasant surprises appear when we look deeper into the
arithmetic of the Legendre curve.  Specifically, we find that the
index of the points we have constructed in the full Mordell-Weil group
is a power of $p$, and the square of this index is the order of a
Tate-Shafarevich group.  This ``class number formula'' is reminiscent
of formulae appearing in the contexts of cyclotomic units and of
Heegner points.  See Theorem~\ref{thm:main} for a
precise statement of our main results.

Here is an outline of the paper: In Sections~\ref{s:E} and
\ref{s:points}, we write down the relevant curve and a collection of
points on it, and in Section~\ref{s:lower-bound} we prove that the
rank of the group they generate is large.  The arguments require only
elementary aspects of the theory of elliptic curves and Galois theory.
Using slightly more advanced ideas, (namely Galois cohomology), in
Section~\ref{s:upper-bound} we then calculate the exact rank of the
Mordell-Weil group and see that our points generate a subgroup of
finite index.  The tools used up to this point will be familiar to
anyone who has studied Silverman's book \cite{SilvermanAEC}.  In
Section~\ref{s:torsion} we compute the torsion subgroup, and in
Sections~\ref{s:Neron} through \ref{s:sha}, we use more sophisticated
techniques (heights, $L$-functions, and the Birch and Swinnerton-Dyer
formula) to bound the index and relate it to the Tate-Shafarevich
group of $E$.  Finally, in Section~\ref{s:Berger} we make a connection
with our paper \cite{Ulmer13a} and Berger's construction, which
yields strong results on the Birch and Swinnerton-Dyer conjecture for
$E$.

There are several avenues for further exploration of the arithmetic of
the Legendre curve, including a more general rank formula, more
explicit points, and a full calculation of the Tate-Shafarevich group.
These are stated more precisely at the end of Section~\ref{s:summary}
and will be taken up in subsequent publications.

This note came about as a result of conversations with Dick Gross at
the 2009 IAS/Park City Mathematics Institute.  It is a pleasure to
thank the organizers of the meeting for the opportunity to give a
course, the audience for their attention and comments, and Gross for
his interest and incisive questions.  Thanks are also due to Chris
Hall for thought-provoking data on $L$-functions and stimulating
conversations on questions related to the topics of this paper, and to
Lisa Berger and Alice Silverberg for several useful comments.

The exposition of the final version of the paper was influenced by
the possibility of generalizations to higher genus.  These were worked
out with an energetic group of mathematicians starting at an AIM
workshop on ``Cohomological methods in abelian varieties'' in March,
2012.  It is a pleasure to thank AIM, the organizers, and the
participants for their collaboration, the fruits of which will appear
in a future paper.

\section{The elliptic curve $E$}\label{s:E}
We refer to \cite{Ulmer11} and the references there for basic results
on elliptic curves over function fields.  Throughout the paper, $p$
will be an odd prime number and $E$ will be the elliptic curve
\begin{equation}\label{eq:E}
E:\quad y^2=x(x+1)(x+t)
\end{equation}
over $\Fp(t)$ or one of its extensions.  The $j$-invariant of $E$ is 
$$j(E)=\frac{256(t^2-t+1)^3}{t^2(t-1)^2}$$
and the discriminant of the model \eqref{eq:E} is
$$\Delta=16t^2(t-1)^2.$$

A word on the $+$ signs in the definition of $E$ is in order.
The term ``Legendre curve'' usually refers to the elliptic curve
$$E_L:\quad y^2=x(x-1)(x-t).$$
The curve $E$ is isomorphic to the quadratic twist 
$$(-1)y^2=x(x-1)(x-t)$$
of $E_L$ (via the change of coordinates $x\mapsto-x$).  Whenever the
ground field has a primitive 4-th root of unity $i$, there is an
isomorphism $E\cong E_L$ which sends $y$ to $iy$.  In most of our
work, the fourth roots of unity will be in the ground field.  Even
when they are not, it turns out that $E$ is more convenient than $E_L$,
and so we use $E$ throughout and refer to it as the ``the Legendre
curve.''

We will consider $E$ and its arithmetic over various extensions of
$K=\Fp(t)$.  Specifically, for each power $q$ of $p$ and each integer
$d$ relatively prime to $p$, we have the extension field
$\Fq(t^{1/d})$.  We will also consider $\Fpbar(t^{1/d})$ where
$\Fpbar$ is an algebraic closure of $\Fp$.  When $d$ is fixed, we
sometimes simplify notation by setting $u=t^{1/d}$.  The most
important extensions for our study will be the fields
$K_d=\Fp(\mu_d,t^{1/d})$ where $\mu_d$ denotes the $d$-th roots of
unity in $\Fpbar$.  Note that if $d$ divides $q-1$, then
$\Fq(t^{1/d})$ is a Galois extension of $K_1=\Fp(t)$ whose Galois
group is a semidirect product of two cyclic groups.  We fix a
primitive $d$-th root of unity in $\Fpbar$ and denote it $\zeta_d$.

Note that as an abstract field, $\Fq(t^{1/d})=\Fq(u)$ is a rational
function field over $\Fq$ and so is isomorphic to $\Fq(t)$.  Thus
instead of fixing our curve $E$ and varying the ground field, we could
fix the ground field $\Fq(t)$ and vary the curve $E$.  The former
point of view will be more convenient for our purposes.

As an aside for experts, very little of what we do requires that $d$
be relatively prime to $p$, but we assume this throughout for
simplicity.

\section{Several points on $E$}\label{s:points}
In this section we write down several points on $E$ rational over
various extensions of $\Fp(t)$.

First we have the points of order 2 visible immediately from the
Weierstrass equation~\eqref{eq:E}: let $Q_0=(0,0)$, $Q_1=(-1,0)$,
and $Q_t=(-t,0)$.

Next assume that $d$ is an integer of the form $d=p^f+1$ with $f$ a
non-negative integer.  Over $K_d=\Fp(\mu_d,u)$ with $u^d=t$, $E$ has the
Weierstrass equation
$$y^2=x(x+1)(x+u^{p^f+1})$$
and the rational point $P(u)=(u,u(u+1)^{d/2})$.  Indeed,
evaluating the right hand side of \eqref{eq:E} at $x=u$, we have
$$u(u+1)(u+u^{p^f+1})=u^2(u+1)(1+u^{p^f}) =u^2(u+1)^{p^f+1}=u^2(u+1)^d.$$

Since $E$ is defined over $K=\Fp(t)\subset K_d$, applying elements of
$\gal(K_d/K)$ to $P(u)$ we get a collection of
points
\begin{equation}\label{eq:points}
  P_i=P(\zeta_d^iu)=(\zeta_d^iu,\zeta_d^iu(\zeta_d^iu+1)^{d/2})\qquad (i\in\Z/d\Z).
\end{equation}

Occasionally we will want to discuss more than one value of $d$ at a
time, and when we do so we label the points above $P_i^{(d)}$ to
emphasize that they are defined over $K_d$.  It is not yet evident
what connections there are between the $P_i^{(d)}$ for varying values of
$d$.
 
We will see by elementary means in the next section that the points
$P_i$ ($i\in\Z/d\Z$) are almost independent---they generate a
subgroup of $E(K_d)$ of rank $d-2$.  Using somewhat more technology,
we will see in the following sections that they generate a subgroup of
finite index, i.e., that the rank of $E(K_d)$ is $d-2$.

One case can be dealt with immediately: Suppose that $f=0$ so
that $d=2$.  Then it is easy to see using the explicit formulas for
addition on $E$ that $P_0^{(2)}$ and $P_1^{(2)}$ have order 4.  In
fact, we find that $2P_0^{(2)}=Q_0$, and four points of order exactly 4
on $E(K_d)$ are:
$$P^{(2)}_0,\quad-P_0^{(2)}=P_0^{(2)}+Q_0,\quad 
P_1^{(2)}=P_0^{(2)}+Q_1,\quad\text{and }-P_1^{(2)}=P_0^{(2)}+Q_t.$$ We
will see later (Prop.~\ref{prop:torsion}) that the torsion in
$E(\Fpbar(t^{1/d}))$ is generated by $Q_0$ and $Q_1$ when $d$ is odd,
and by $P^{(2)}_0$ and $Q_1$ when $d$ is even.

Using the points $P^{(2)}_i$, we can see by elementary means certain
relations among the $P_i^{(d)}$.
\begin{prop}\label{prop:relations}
  Suppose that $f\ge1$ and $d=p^f+1$.  Then we have the following
  equalities in the Mordell-Weil group $E(K_d)$:
$$\sum_{i=0}^{d-1}P_i^{(d)}=Q_t
\quad\text{and}\quad\sum_{i=0}^{d-1}(-1)^iP_i^{(d)}=Q_1.$$ 
\end{prop}

\begin{proof}
  Consider the rational function $f=y-x(x+1)^{d/2}$ on $E$.  Because
  $f$ is a polynomial in $x$ and $y$, it is regular away from $O$, the
  origin of $E$, and it has a pole there of order $d+2$.  It is easy
  to see that $f$ vanishes at each $P_i$ ($i\in\Z/d\Z$) and at $Q_0$
  and $Q_1$.  (Here and below, we write $P_i$ for $P^{(d)}_i$.)  Since
  the degree of the divisor of $f$ is zero, $f$ must vanish simply at
  each of these points.  Thus we have
$$\dvsr(f)=\sum_{i\in\Z/d\Z}[P_i]+[Q_0]+[Q_1]-(d+2)[O].$$
Abel's theorem then yields
$$\sum_{i\in\Z/d\Z}P_i=-Q_0-Q_1=Q_t$$
in $E(K_d)$.

The second equality in the statement of the proposition follows
similarly from a consideration of the rational function
$g=x^{d/2-1}y-t^{1/2}(x+1)^{d/2}$ whose divisor is
$$\dvsr(g)=\sum_{i\in\Z/d\Z}[(-1)^iP_i]+[Q_1]-(d+1)[O].$$  
We leave the details of this case as an exercise for the reader.
\end{proof}

\begin{rem}
Adding and subtracting the two relations above,
we find that 
$$2\sum_{i\in2\Z/d\Z}P_i=2\sum_{i\in2\Z/d\Z}P_{i+1}=Q_0.$$
Thus the two sums are points of order 4.  We will see below
(Prop.~\ref{prop:trace-to-d=2}) that
$$\sum_{i\in2\Z/d\Z}P_i=P^{(2)}_{1+d/2}\quad\text{and}\quad
\sum_{i\in2\Z/d\Z}P_{i+1}=P^{(2)}_{d/2}.$$
These will turn out to be the only relations among the $P_i$ (for a
fixed value of $d$).
\end{rem}

\section{A lower bound on the rank}\label{s:lower-bound}
In this section we give a very simple proof that the rank of $E(K_d)$
is unbounded as $d$ varies.  Experts will see that the argument uses a
coboundary map in Galois cohomology, but to make this section as
accessible as possible we argue from first principles.

\begin{prop}\label{prop:coboundary}
  Let $d$ be a positive integer not divisible by $p$ and write $\zeta_d$
  for a fixed primitive $d$-th root of unity in $\Fpbar$. Let
  $F=\Fq(t^{1/d})=\Fq(u)$. Then the map
  $E(F)\to (\Z/2\Z)^d$ which sends
$$P=(x,y)\mapsto\begin{cases}
\left(\ord_{u=\zeta_d^{j}}(x+1)\pmod2\right)_{j=0,\dots,d-1}&\text{if
  $P\neq(-1,0)$ and $P\neq O$}\\
\left(1,\dots,1\right) &\text{if $P=(-1,0)$}\\
\left(0,\dots,0\right) &\text{if $P=O$}
\end{cases}
$$
induces a homomorphism $\nu:E(F)/2E(F)\to(\Z/2\Z)^d$.  Suppose
further that $d=p^f+1$ and let $q$ be a power of $p$ such that
$d$ divides $q-1$. Then the homomorphism $\nu$ maps the subgroup
of $E(F)$ generated by the points $P_i$
\textup{(}$i\in\Z/d\Z$\textup{)} surjectively onto $(\Z/2\Z)^d$.
\end{prop}

\begin{proof}
Suppose we are given an elliptic curve $E'$ in Weierstrass form
$$Y^2=X^3+a_2X^2+a_4X$$
over a field $k$ of characteristic $\neq2$.  Note that the point
$(0,0)$ on $E'$ is of order 2.  We claim that the map $E'(k)\to
k^\times$ that sends
$$P=(X,Y)\mapsto\begin{cases}
X&\text{if $X\neq0$ and $P\neq O$}\\
a_4&\text{if $X=0$}\\
1&\text{if $P=O$}
\end{cases}
$$
induces a homomorphism $E'(k)/2E'(k)\to k^\times/k^{\times 2}$ (where
$k^\times/k^{\times 2}$ denotes the multiplicative group of $k$ modulo
squares).  Indeed, if we have three points in $E'(k)$ with
$P_1+P_2=P_3$, then $P_1$, $P_2$, and $-P_3$ lie on a line $Y=\lambda
X+\nu$, and the $X$-coordinates $X_i=X(P_i)$ are the three roots of
$$(\lambda X+\nu)^2=X^3+a_2X^2+a_4X$$
or equivalently of
$$0=X^3+(a_2-\lambda^2)X^2+(a_4-2\lambda\nu)X-\nu^2.$$
If none of the $X_i$ are 0, it follows that $X_1X_2X_3=\nu^2$ and so
$X_3=X_1X_2$ in $k^\times/k^{\times 2}$.  If one of the $X_i=0$, say 
$X_1=0$ for concreteness, then we have $X_2X_3=a_4$ and so
$X_3=a_4X_2$ in $k^\times/k^{\times 2}$.  The cases where $X_2=0$ or
$X_3=0$ are similar.  This shows the we have a homomorphism $E'(k)\to
k^\times/k^{\times 2}$.  Since the target is a group of exponent 2, we
get an induced homomorphism $E'(k)/2E'(k)\to k^\times/k^{\times 2}$ as
claimed.

Now we apply this observation to the Legendre curve over $k=F$ with the
change of coordinates $X=x+1$, $Y=y$.  We get a homomorphism
$E(F)/2E(F)\to F^\times/F^{\times2}$ defined by
$$P=(x,y)\mapsto\begin{cases}
x+1&\text{if $x\neq-1$}\\
1-t&\text{if $x=-1$}
\end{cases}
$$

Next we compose with the homomorphism
$F^\times/F^{\times2}\to(\Z/2\Z)^d$ that sends the class of $z\in
F^\times$ to $\left(\ord_{u=\zeta_d^{j}}(z)\pmod
  2\right)_{j=0,\dots,d-1}$. Noting that $1-t=1-u^d$ vanishes simply
at each of the points $u=\zeta_d^{j}$, we see that we have a
homomorphism $\nu:E(F)/2E(F)\to(\Z/2\Z)^d$ given by the formula in
the statement of the proposition.

Now assume that $d=p^f+1$ and $d$ divides $q-1$ so that the
points $P_i$ lie in $E(F)$.  Then it is easy to see that that the point
$P_i=(\zeta_d^iu,\zeta_d^iu(\zeta_d^iu+1)^{d/2})$ maps to the vector
in $(\Z/2\Z)^d$ with a 1 in the position labeled by $j$ with
$j\equiv(d/2)-i\pmod d$ and zeroes elsewhere.  Thus the subgroup of
$E(F)$ generated by the $P_i$ maps surjectively onto $(\Z/2\Z)^d$.

This completes the proof of the proposition.
\end{proof}

For the statement of the next result, suppose that $d=p^f+1$ and let
$G$ be the quotient of the free abelian group $\Z^d$ by the subgroup
generated by the vectors
$$(2,-2,2,-2\dots,2,-2)$$ 
(with 2's and $-2$'s alternating) and 
$$(2,2,\dots,2)$$ 
(all entries equal to 2).  We have
\begin{equation}\label{eq:Gstructure}
G\cong\Z^{d-2}\oplus\Z/2\Z\oplus\Z/4\Z.
\end{equation}
By Proposition~\ref{prop:relations}, the map $\Z^d\to E(K_d)$,
$(a_0,\dots,a_{d-1})\mapsto\sum_ia_iP_i$ factors through $G$.
Moreover, the 2-torsion points of $E(K_d)$ are all images of
2-torsion elements of $G$.

\begin{cor}\label{cor:lower-bound}
  For $d=p^f+1$ and $K_d=\Fp(\mu_d,t^{1/d})$, we have $\rk E(K_d)\ge
  d-2$. More precisely, the subgroup of $E(K_d)$ generated by the
  $P_i$ \textup{(}$i\in\Z/d\Z$\textup{)} is isomorphic to $G$.
\end{cor}

\begin{proof} 
  By definition, the subgroup of $E(K_d)$ generated by the $P_i$ is
  the image of the homomorphism $G\to E(K_d)$. On the other hand,
  since $d=p^f+1$ and $K_d$ contains the $d$-th roots of unity,
  Proposition~\ref{prop:coboundary} shows that the composed map
  $G/2G\to E(K_d)/2E(K_d)\labeledto{\nu}(\Z/2\Z)^d$ is surjective, and
  therefore an isomorphism. Thus $G/2G\to E(K_d)/2E(K_d)$ is
  injective. If $a=(a_0,\dots,a_{d-1})\in\Z^d$ maps to zero in
  $E(K_d)$, i.e., if $\sum_ia_iP_i=O$, then it follows that $a$ is
  zero in $G/2G$, i.e., $a=2a'$ for some element $a'\in G$.  Now $a'$
  maps to a 2-torsion point in $E(K_d)$, and since the 2-torsion of
  $G$ maps onto the 2-torsion of $E(K_d)$, we may modify $a'$ so that
  we still have $a=2a'$ and $a'$ goes to zero in $E(K_d)$.  But now
  $a'=2a^{\prime\prime}$ for some $a^{\prime\prime}$ in $G$.
  Iterating the argument, we see that $a$ is divisible in $G$ by
  arbitrarily large powers of 2.  But we see from
  equation~\eqref{eq:Gstructure} that $G$ is finitely generated, so
  has no non-trivial divisible elements.  This implies that $a=0$.
  Thus $G\to E(K_d)$ is injective and this finishes the proof of the
  corollary.
\end{proof}

Using Galois theory, we can estimate the rank of $E(\Fq(t^{1/d}))$ for any
$q$ and any divisor $d$ of $p^f+1$.

\begin{cor}\label{cor:more-ranks}
If $d$ divides $p^f+1$ for some $f$, $q$ is any power
of $p$, and $F=\Fq(t^{1/d})$, we have
$$\rk E(F)\ge
\sum_{\substack{e|d\\e>2}}\frac{\varphi(e)}{o_q(e)}$$
where the sum is over divisors $e$ of $d$ greater than 2, $\varphi$ is
Euler's function, and $o_q(e)$ is the order of $q$ in
$(\Z/e\Z)^\times$.  
\end{cor}

\begin{rem}
  The Corollary gives an elementary and constructive proof that for
  every odd $p$ there is a non-isotrivial elliptic curve over the
  rational function field $\Fp(u)$ with arbitrarily large rank.
  Indeed, if we take $d=p^f+1$ and $q=p$, then $o_q(e)\le2f$ for all
  divisors $e$ of $d$.  Thus over $F=\Fp(t^{1/d})\cong\Fp(u)$, we have
  $\rk E(F)\ge(p^f-1)/(2f)$.  The same result is also known for $p=2$;
  see \cite[8.1 and 8.2]{Ulmer13a}.
\end{rem}

\begin{proof}[Proof of Corollary~\ref{cor:more-ranks}]
  Choose an $f$ such that $d$ divides $d'=p^f+1$ and $q$ divides
  $|\Fp(\mu_{d'})|=p^{2f}$.  Then $F$ is a subfield of $K_{d'}=\Fp(\mu_{d'},t^{1/d'})$
  and $K_{d'}/F$ is a Galois extension with Galois group
  $G=\gal(K_{d'}/F)$.  We have rational points $P_i\in E(K_{d'})$ for
  $i\in\Z/d'\Z$.  We will take the trace of these points down to
  $E(F)$ and see that they give the desired lower bound on the rank.
  More precisely we will compute the rank of the subgroup of
  $E(F)=E(K_{d'})^G$ generated by certain Galois invariant elements of
  the subgroup of $E(K_{d'})$ generated by the $P_i$.

  The Galois group of $K_{d'}$ over $F$ is a semidirect product.  More
  precisely, choosing a primitive $d'$-th root of unity $\zeta$, we
  have an identification:
$$\gal(K_{d'}/F)\cong(d\Z/d'\Z)\sdp \langle\Fr_q\rangle.$$
Here $\langle\Fr_q\rangle$ is the cyclic group of order
$[\F_q(\mu_{d'}):\Fq]$ generated by the $q$-power Frobenius acting
trivially on $t^{1/d'}$ and as usual on $\Fq(\mu_{d'})\subset K_{d'}$, and
$jd\in d\Z/d'\Z$ acts trivially on $\Fq(\mu_{d'})\subset K_{d'}$ and sends
$t^{1/d'}$ to $\zeta^{jd}t^{1/d'}$.

Looking back at the definition of $P_i$, it is clear the orbit of $P_i$
under $\gal(K_{d'}/F)$ consists of all points $P_j$ where $j\equiv
q^ai\pmod d$ for some power $q^a$ of $q$.  Summing all the points
in the orbit yields a rational point in $E(F)$.  If $o\subset\Z/d'\Z$
is an orbit, we write $P_o$ for the corresponding sum.

To estimate the rank, we first note that the number of orbits is
$$r=\sum_{{e|d}}\frac{\varphi(e)}{o_q(e)}.$$  
(The orbit through $P_i$ is counted in the term corresponding to
$e=d/\gcd(d,i)$.)  It is easy to see, using
Proposition~\ref{prop:coboundary}, that the points just constructed
are linearly independent in $E(K)/2E(K)\cong(\Z/2\Z)^{d'}$.  Thus the
rank of the subgroup of $E(F)$ they generate is $r-s$ where $2^s$ is
the number of 2-torsion points in the span of the $P_o$.

By Proposition~\ref{prop:relations}, the 2-torsion point $Q_t$ is the
sum over all orbits $o$ of $P_o$.  Again using
Proposition~\ref{prop:relations}, it is easy to see that the 2-torsion
point $Q_0$ is in the span of the $P_o$ if and only if $d$ is even.
Thus $s=1$ if $d$ is odd and $s=2$ if $d$ is even.  This gives the
formula for the rank of $E(F)$ in the statement of the corollary.
\end{proof}

\begin{exs}\mbox{}
\begin{enumerate}
\item Suppose that $p\equiv3\pmod4$,  $d=4$, and $q=p$.  Let
  $d'=p+1$, $F=\Fp(t^{1/4})$, and $K=K_{d'}$.  We have points $P_i\in
  E(K)$ ($i\in\Z/d'\Z$).  There are three orbits of $\gal(K/F)$ on the
  $P_i$: 
\begin{align*}
o_0&=\{P_i|i\in\Z/d'\Z\text{ and }i\equiv0\pmod4\}\\
o_{13}&=\{P_i|i\in\Z/d'\Z\text{ and }i\equiv1\pmod4\text{ or }i\equiv3\pmod4\}\\
o_2&=\{P_i|i\in\Z/d'\Z\text{ and }i\equiv2\pmod4\}
\end{align*}
We have $P_{o_0}+P_{o_{13}}+P_{o_2}=Q_t$ and
$P_{o_0}-P_{o_{13}}+P_{o_2}=Q_0$, and the subgroup of $E(F)$ generated by
$P_{o_0}$, $P_{o_{13}}$, and $P_{o_2}$ has rank 1.
\item Now suppose that $p\equiv3\pmod4$, $d=4$, and $q=p^2$.  The
  situation is similar to the previous example, except that the orbit
  $o_{13}$ breaks up into two orbits
\begin{align*}
o_1&=\{P_i|i\in\Z/d'\Z\text{ and }i\equiv1\pmod4\}\\
o_3&=\{P_i|i\in\Z/d'\Z\text{ and }i\equiv3\pmod4\}
\end{align*}
and the subgroup of $E(F)$ generated by $P_{o_0},\dots,P_{o_3}$ has
rank 2.
\item Suppose now that $d=d'=p+1$ and $q=p$.  Let
  $K=K_{d'}=\F_{p^2}(t^{1/(p+1)})$ and $F=\Fp(t^{1/(p+1)})$.  Then
  $\gal(K/F)$ has order 2, generated by the Frobenius $\Fr_p$.  The
  orbits of Galois on the $P_i$ are $\{P_0\}$,
  $\{P_{d/2}\}$, and $\{P_i,P_{-i}\}$ for $1\le i\le(p-1)/2$.  We
  find that the subgroup of $E(F)$ generated by the $P_o$ has rank
  $(p-1)/2$.
\item We can make the previous example more explicit.  For an integer
  $i$ in the range $1\le i\le (p-1)/2$, let
  $b=\zeta_d^i+\zeta_d^{pi}=\zeta_d^i+\zeta_d^{-i}\in\Fp$ so that the
  roots of $x^2-bx+1$ are $\zeta_d^i$ and $\zeta_d^{pi}$.  Then
  computing with the standard formulae leads to an expression for
  $P_{\{i,-i\}}=P_i+P_{-i}$ in terms of $b$.  More precisely, we find
  that the $x$ coordinate of $P_{\{i,-i\}}$ is
$$\frac{2u^{p+1}+bu^p+bu+2-2(u^2+ub+1)^{d/2}}{b^2-4}.$$
(We omit the $y$ coordinate, which is even more unwieldy.)  This
example suggests that the coordinates of the points $P_o$ may in
general be rather complicated.
\end{enumerate}
\end{exs}

\section{An upper bound on the rank}\label{s:upper-bound}
In this section we compute the Selmer group of $E$ for multiplication
by 2 and complete the calculation of the rank of $E$ when $d$ divides
$p^f+1$ for some $f$.

Throughout, we fix a value of $d$ (prime to $p$ as always) and a power
$q$ of $p$.  Let $F=\Fq(t^{1/d})=\Fq(u)$ and let $\Fsep$ be a
separable closure of $F$.  Write $G_F$ for $\gal(\Fsep/F)$.

We have an exact sequence of $G_F$-modules:
\begin{equation}\label{eq:2-isog}
0\to E[2]\to E(\Fsep)\labeledto{2}E(\Fsep)\to0
\end{equation}
where ``$2$'' denotes the multiplication-by-2 homomorphism and
$E[2]\subset E(F)\subset E(\Fsep)$ is its kernel.

Taking Galois cohomology leads to an injection
$$0\to E(F)/2E(F)\labeledto{\lambda} H^1(G_F,E[2]).$$

For each place $v$ of $F$, we have the completion $F_v$, its separable
closure $\Fsep_v$, and the Galois group $G_{F_v}=\gal(\Fsep_v/F_v)$.
There is a short exact sequence analogous to \eqref{eq:2-isog} and an
injection
$$0\to E(F_v)/2E(F_v)\labeledto{\lambda_v} H^1(G_{F_v},E[2]).$$

The field inclusion $F\subset F_v$ and the induced restriction maps of
Galois groups $G_{F_v}\into G_F$ lead to a diagram
\begin{equation*}
\xymatrix{
0\ar[r]&E(F)/2E(F)\ar[r]^{\lambda}\ar[d]&H^1(G_F,E[2])\ar[d]^r\\
0\ar[r]&\prod_v
E(F_v)/2E(F_v)\ar[r]^{\prod_v\lambda_v}&\prod_vH^1(G_{F_v},E[2]).}
\end{equation*}

We define the \emph{$2$-Selmer group} $\Sel_2(E/F)$ to be
$r^{-1}(\im(\prod_v\lambda_v)$.  In other words, $\Sel_2(E/F)$ consist
of classes $c\in H^1(G_F,E[2])$ such that for all $v$, the
component of $r(c)$ indexed by $v$ lies in the image of $E(F_v)/2E(F_v)\to
H^1(G_{F_v},E[2])$.  Since the diagram above is commutative, we have
an injection
$$0\to E(F)/2E(F)\labeledto{\lambda}\Sel_2(E/F).$$

Choosing a basis of $E[2]$ yields an identification
$E[2]\cong\mu_2^2$ of $G_F$-modules.  In the calculations below, we
always use the basis $\{Q_0,Q_1\}$.  Kummer theory yields an
isomorphism $H^1(G_F,\mu_2)\cong F^\times/F^{\times2}$ and so the
Selmer group $\Sel_2(E/F)$ is a subgroup of
$(F^\times/F^{\times2})^2$.

We define a homomorphism 
$$\lambda':\Sel_2(E/F)\to(\Z/2\Z)^d$$
as the composition of the inclusion $\Sel_2(E/F)\subset(F^\times/F^{\times2})^2$
and the map
\begin{align*}
(F^\times/F^{\times2})^2&\to (\Z/2\Z)^d\\
[a,b]&\mapsto\left(\ord_{u=\zeta_d^j}b\pmod2\right)_{j=0,\dots,d-1}.
\end{align*}

\begin{prop}\label{prop:Sel}
Recall that $F=\Fq(t^{1/d})=\Fq(u)$ and $E$ is the elliptic curve
  given by equation~\eqref{eq:E}.  Assume that $d$ divides $q-1$.
\begin{enumerate}
\item The composed homomorphism 
$$E(F)/2E(F)\labeledto{\lambda}\Sel_2(E/F)\labeledto{\lambda'}(\Z/2\Z)^d$$ 
is equal to the map $\nu$ appearing in
Proposition~\ref{prop:coboundary}.
\item If $d$ is even then $\lambda'$ is injective.  Its image has
  order $2^d$ if $2d$ divides $q-1$, order $2^{d-2}$ if $(q-1)/d$ is
  odd and $q\equiv1\pmod 4$, and order $4$ if $q\equiv3\pmod4$.
\item If $d$ is odd and $q\equiv1\pmod4$ then $\lambda'$ is surjective
  with kernel of order $2$.  If $d$ is odd and $q\equiv3\pmod4$ then
  $\lambda'$ is injective with image of order $4$.
\end{enumerate}
\end{prop}

Before proving the proposition, we state several corollaries.

\begin{cor}\label{cor:rks} \mbox{}
\begin{enumerate}
\item If $d$ is even, then $\rk E(F)\le d-2$.
\item If $d$ is odd, then $\rk E(F)\le d-1$.
\item If $d=p^f+1$ and $d$ divides $q-1$, then $\rk E(F)=d-2$.
\item More generally, if $d$ divides $p^f+1$ for some $f$, then 
$$\rk E(F)=
\sum_{\substack{e|d\\e>2}}\frac{\varphi(e)}{o_q(e)}$$
where the sum is over divisors $e$ of $d$ greater than 2, $\varphi$ is
Euler's function, and $o_q(e)$ is the order of $q$ in
$(\Z/e\Z)^\times$.
\end{enumerate}
\end{cor}

\begin{proof}[Proof of Corollary~\ref{cor:rks}]
  For the upper bounds, we may extend the finite field $q$ and so we
  may assume that $2d$ divides $q-1$ when $d$ is even and $4d$
  divides $q-1$ when $d$ is odd.  Proposition~\ref{prop:Sel}
  then shows that $\Sel_2(E/F)$ has cardinality $2^d$ when $d$ is even
  and $2^{d+1}$ when $d$ is odd.  Since $E$ has two independent points
  of order $2$ defined over $\Fp(t)\subset F$, it follows that $\rk
  E(F)$ is bounded above by $d-2$ when $d$ is even and by $d-1$ when
  $d$ is odd.  This establishes the first two claims.  The third claim
  is immediate from the first and the lower bound in
  Corollary~\ref{cor:lower-bound}.

  To establish the last claim, choose a value of $f$ so that $d$
  divides $d'=p^f+1$ and $q$ divides $p^{2f}$, in which case $F$ is a
  subfield of $K_{d'}=\Fp(\mu_{d'},t^{1/d'})$.  Then $K_{d'}$ is
  Galois over $F$ and we write $G$ for $\gal(K_{d'}/F)$.  Write
  $V_{d'}$ for the subgroup of $E(K_{d'})$ generated by $P^{(d')}_i$
  ($i\in\Z/d'\Z$).  The lower bound of Corollary~\ref{cor:lower-bound}
  together with Proposition~\ref{prop:Sel} shows that $V_{d'}$ has
  finite index in $E(K_{d'})$.  Taking Galois invariants, we find that
  $V^{G}_{d'}$ has finite index in $E(F)=E(K_{d'})^G$.  On the other
  hand, the proof of the lower bound of Corollary~\ref{cor:more-ranks}
  constructs a submodule of $V_{d'}^G$ of finite index and rank equal
  to $\sum_{\substack{e|d\\e>2}}\frac{\varphi(e)}{o_q(e)}$.  This
  gives the rank asserted in the Corollary.
\end{proof}

\begin{proof}[Proof of Proposition~\ref{prop:Sel}]
  We use \cite{BrumerKramer77} as a general reference for 2-descent;
  it has fairly complete results on the 2-Selmer group for elliptic
  curves over local and global fields of characteristic $\neq2$.

  To prove the first part of the proposition, we note from
  \cite[\S2]{BrumerKramer77} that under the identification
  $H^1(F,E[2])\cong (F^\times/F^{\times2})^2$, the coboundary
  homomorphism $E(F)\to H^1(F,E[2])$ becomes the map
\begin{align*}
E(F)&\to (F^\times/F^{\times2})^2\\
P=(x,y)&\mapsto\begin{cases}
[x,x+1]&\text{if $P\neq Q_0$, $Q_1$, $O$}\\
[t,1]&\text{if $P=Q_0$}\\
[-1,1-t]&\text{if $P=Q_1$}\\
[1,1]&\text{if $P=O$}
\end{cases}
\end{align*}
(Here and below we use square brackets to denote the class of an
element $a\in F^\times$ or a pair $(a,b)\in(F^\times)^2$ in a quotient
group.)  It is immediate from this formula that the composition
$\lambda'\lambda$ is equal to $\nu$.

For the second and third parts of the proposition, we need to compute
the Selmer group explicitly.  For each place $v$ of $F$, we write
$\Sel_2(E/F_v)$ for the local Selmer group, i.e., for the image of the
coboundary map 
$$E(F_v)\to H^1(F_v,E[2])\cong (F_v^\times/ F_v^{\times2})^2.$$
Recall that the ambient group $(F_v^\times/ F_v^{\times2})^2$ has
order 16.  By \cite[Lemma~3.1]{BrumerKramer77}, the local Selmer
group has order 4.  

In order to apply the local results of \cite{BrumerKramer77}, we need
to know the reduction types of $E$ at the places of $F$.  Using Tate's
algorithm \cite{Tate75}, we find that $E$ has good reduction away from
$u=0$, $u\in\mu_d$, and $u=\infty$.  At $u=0$ the reduction is split
multiplicative for any value of $q$.  At places with $u\in\mu_d$,
the reduction is split multiplicative when $q\equiv1\pmod4$, and
non-split multiplicative when $q\equiv3\pmod4$.  At $u=\infty$, the
reduction is additive when $d$ is odd and split multiplicative (for
any $q$) when $d$ is even.

To describe $\Sel_2(E/F_v)$ explicitly, let $\O_v$ be the ring of
integers of $F_v$ and let $\pi_v$ be a uniformiser at $v$.  At places
where $E$ has good reduction, \cite[Cor.~3.3]{BrumerKramer77}
says that the local Selmer group is
$$\Sel_2(E/F_v)=(\O_v^\times/\O_v^{\times2})^2\subset (F_v^\times/
F_v^{\times2})^2.$$

The local Selmer group at places of split multiplicative reduction is
given by \cite[Prop.~4.1]{BrumerKramer77}.  That result tells us that
at $u=0$ the local group is $(F_v^\times/ F_v^{\times2})\times\{1\}$,
at places with $u\in\mu_d$, when $q\equiv1\pmod4$ the local group is
$\{1\}\times(F_v^\times/ F_v^{\times2})$, and at $u=\infty$, when $d$
is even the local group is
$$\left\{[a,b]\in (F_v^\times/ F_v^{\times2})^2\left|\,a=b\right.\right\}.$$

The local Selmer group at places of non-split multiplicative reduction is
given by \cite[Prop.~4.3]{BrumerKramer77}.  That result tells us that
at places with $u\in\mu_d$, when $q\equiv3\pmod4$ the local group is
is generated by $[-1,\pi_v]$ and $[-1,-\pi_v]$ where $\pi_v$ is a
uniformiser at $v$.

Unfortunately, \cite{BrumerKramer77} does not give the local Selmer
group in general for places of additive reduction.  However, we can
determine it by exhibiting explicit classes.  Namely, if we apply the
coboundary map
$$E(F)\to (F^\times/F^{\times2})^2\to (F_\infty^\times/F_\infty^{\times2})^2$$
at the place $u=\infty$ to the points $Q_0$ and $Q_1$, we obtain
non-trivial elements $[t,1]=[\pi_\infty,1]$ and
$[-1,1-t]=[-1,-\pi_\infty]$.  (Here the square brackets indicates
classes in $(F_\infty^\times/F_\infty^{\times2})^2$.)  Since the local Selmer
group has order 4, it must be generated by these classes.  
Note that if $q\equiv1\pmod4$ then $-1$ is a square in $\Fq$ and
$[-1,-\pi_\infty]=[1,\pi_\infty]$.

We can now describe the global Selmer group when $d$ is even.  First
assume that $2d$ divides $q-1$.  We use
crucially that the Picard group of $\P^1$ is trivial, so that every
divisor of degree 0 is the divisor of a function.  If $[a,b]$ is in the
Selmer group, then the divisor of $b$ has the form $D_1+2D_2$ where
$D_1$ is supported on the places $u\in\mu_d$ and $u=\infty$.
Modifying $b$ by a square, we may assume that $b$ is a polynomial in
$u$ whose zeroes are all in $\mu_d$ and that these zeroes have order 0
or 1.  The local condition at $u=0$
then requires that the value of $b$ at $u=0$ should be a square;
modifying by the square of an element of $\Fq^\times$, we may assume
that the value of $b$ at $u=0$ is 1.  We have thus put $b$ in the
standard form
$$b=\prod_{j=0,\dots, d-1}(1-u\zeta_d^j)^{e_j}$$
where $(e_0,\dots,e_{d-1})$ is an arbitrary tuple in $\{0,1\}^d$.  Now
consider the local conditions on $a$.  Modifying $a$ by a square, we
may assume that $a$ is regular away from $0$ and $\infty$, so it has
the form $cu^f$ where $c\in\Fq$ and where we may take $f=0$ or $f=1$.
Since the value of $a$ must be a square at places where $u\in\mu_d$,
and since our hypotheses imply that a $d$-th root $\zeta_d$ is a
square in $\Fq$, we must have that $c$ is a square; modifying $a$ by
the square of a scalar, we may take $c$ to be 1.  Finally, the local
condition at $\infty$ implies that $f\equiv\sum_je_j\pmod 2$.
Summarizing, we have that when $d$ is even and $2d$ divides $q-1$,
$\Sel_2(E/F)$ consists of classes $[a,b]\in (F^\times/F^{\times2})^2$
represented by pairs $(a,b)$ of the form
$$a=u^f,\quad b=\prod_{j=0,\dots, d-1}(1-u\zeta_d^j)^{e_j}$$
where $(e_0,\dots,e_{d-1})\in\{0,1\}^d$ is arbitrary and $f\in\{0,1\}$
satisfies $f\equiv\sum_je_j\pmod2$.  It is immediate that the map
\begin{align*}
\Sel_2(E/F)&\labeledto{\lambda'}(\Z/2Z)^d\\
[a,b]&\mapsto\left(\ord_{u=\zeta_d^j}b\pmod2\right)_{j=0,\dots,d-1}\
\end{align*}
is an isomorphism.  

The argument when $d$ is even, $(q-1)/d$ is odd, and $q\equiv1\pmod4$
is very similar, except in this case $\zeta_d$ is not a square in
$\Fq$.  We find that $f$ must be even and $c$ must be a square above
(i.e., we may take $a=1$).  The local conditions at $u=\infty$ then imply
that the allowable tuples $(e_0,\dots,e_{d-1})$ are those that satisfy
$\sum_je_j\equiv\sum_jje_j\equiv0\pmod 2$.  Thus in this case
$\lambda'$ is injective and its image has order $2^{d-2}$.

Now we assume that $q\equiv3\pmod4$.  The argument is the same whether
$d$ is even or odd.  Considering the projection of the local Selmer
group to the first factor $F_v^\times/F_v^{\times2}$, we see that $a$
must be equivalent mod squares to $cu^f$ where $c\in\{\pm1\}$ and
$f\in\{0,1\}$.  Considering the projection to the second factor, we
see that we may again take $b$ of the form
$\prod_j(1-u\zeta_d^j)^{e_j}$.  On the other hand, the local
conditions at $u\in\mu_d$, where we have non-split multiplicative
reduction, link $a$ and $b$.  More precisely, they require that
$\ord_v(b)$ be odd exactly when the value of $a$ at $v$ is a
non-square.  Thus $a$ determines $b$.  It follows that the order of
the global Selmer group is at most 4.  One may argue directly that it
is exactly 4, or note that the images of $Q_0$ and $Q_1$ already
generate a subgroup of order 4.  We note that the image of $\lambda'$
in this case is the subgroup of $(\Z/2\Z)^d$ generated by
$(1,0,1,0,\dots,1,0)$ and $(0,1,0,1,\dots,0,1)$, and $\lambda'$ is
again injective.

The last case to be dealt with is when $d$ is odd and
$q\equiv1\pmod4$.  This case is very similar to the case where $d$ is
even and $2d$ divides $q-1$, the only difference being the local
conditions at $u=\infty$.  By an argument similar to that above, we
see that a class $[a,b]$ in $\Sel_2(E/F)$ can be represented by a pair
$(a,b)$ with
$$a=u^f,\quad b=\prod_{j=0,\dots,d-1}(1-u\zeta_d^j)^{e_j}$$
where $(e_0,\dots,e_{d-1})\in\{0,1\}^d$ and $f\in\{0,1\}$ are arbitrary.
It is immediate that $\lambda'$ is surjective with kernel of order 2.

This completes the proof of the Proposition.
\end{proof}

\section{Torsion}\label{s:torsion}
To complete the elementary part of our analysis of the arithmetic of
$E$, we compute the torsion subgroup of $E(\Fq(t^{1/d}))$.

\begin{prop}\label{prop:torsion}
Let $d$ be a positive integer not divisible by $p$.
If $d$ is odd, then 
$$E(\Fp(t^{1/d}) )_{tor}=E(\Fpbar(t^{1/d}))_{tor}\cong(\Z/2\Z)\times(\Z/2\Z).$$
If $d$ is even, then 
$$E(\Fp(t^{1/d}))_{tor}=E(\Fpbar(t^{1/d}))_{tor}\cong(\Z/2\Z)\times(\Z/4\Z).$$
\end{prop}

\begin{proof}
We have already seen in Section~\ref{s:points} that the torsion
subgroup of $E(\Fp(t^{1/d}))$ is at least as large as what is asserted
in the proposition.  Thus it will suffice to prove that the torsion
subgroup of $E(\Fpbar(t^{1/d}))$ is no larger than what is asserted.

We work one prime at a time, starting with $2$-power torsion.  If $d$
is odd, then the subgroup of $E(\Fpbar(t^{1/d}))$ generated by the
2-torsion points $Q_0$ and $Q_1$ maps injectively to $(\Z/2\Z)^d$
under the map in the display in Proposition~\ref{prop:coboundary} (for
any value of $q$).  Thus no non-trivial combination of these points is
divisible by $2$ in $E(\Fpbar(t^{1/d}))$ and so the 2-power-torsion
subgroup of $E(\Fpbar(t^{1/d}))$ is generated by $Q_0$ and $Q_1$.
Similarly, when $d$ is even, we have a subgroup of
$E(\Fpbar(t^{1/d}))$ isomorphic to $\Z/2\Z\times\Z/4\Z$ generated by
$Q_1$ and $P^{(2)}_0$, call it $H$.  Under the displayed map in
Proposition~\ref{prop:coboundary}, $H/2H$ maps injectively to
$(\Z/2\Z)^d$ (for any value of $q$), so again we find that the 2-power
torsion of $E(\Fpbar(t^{1/d}))$ is generated by $Q_1$ and $P^{(2)}_0$.

Since the $j$-invariant of $E$ is visibly not a $p$-th power in
$\Fpbar(t^{1/d})$, Proposition~I.7.3 of \cite{Ulmer11} implies that
$E$ has no non-trivial $p$-torsion points over that field.

To finish, we must show that for all primes $\ell$ not dividing $2p$,
$E$ has no non-trivial $\ell$-torsion points defined over
$\Fpbar(t^{1/d})$ for any $d$.  For this we may assume that $d$ is
even, and we assume so for the rest of the proof.  Suppose then that
we had a point of order $\ell$.  Then we would have an injection
$\Z/2\Z\times\Z/4\ell\Z\into E(\Fpbar(t^{1/d}))$.  This would lead to
a morphism $\P^1\to X(2,4\ell)$ where $\P^1$ is the curve over
$\Fpbar$ whose function field is $\Fpbar(t^{1/d})\cong\Fpbar(u)$ and
$X(2,4\ell)$ is the modular curve over $\Fpbar$ parameterizing
elliptic curves equipped with an injection
$\Z/2\Z\times\Z/4\ell\Z\into E$.  Since the $j$-invariant of $E$ is
non-constant, the morphism $\P^1\to X(2,4\ell)$ would also be
non-constant.  But that is impossible by the lemma below.  Thus
$E(\Fpbar(t^{1/d}))$ has no non-trivial $\ell$-torsion.
\end{proof}

\begin{lemma}
  If $\ell$ is a prime number not dividing $2p$, then the modular
  curve $X(2,4\ell)$ over $\Fpbar$ parameterizing elliptic curves
  equipped with a subgroup isomorphic to $\Z/2\Z\times\Z/4\ell\Z$ has
  positive genus.
\end{lemma}

\begin{proof}
  The curve $X(2,4\ell)$ covers $X(1,4\ell)=X_1(4\ell)$ and it is well
  known that $X_1(n)$ has positive genus for $n\ge13$ or $n=11$.
  (See, for example, \cite[1.6]{ShimuraIATAF} or \cite[4.2]{MiyakeMF}
  for formulas for the genus of a modular curve and the ramification
  assertions just below.)  This suffices to prove the lemma for
  $\ell>3$.  For $\ell=3$ we consider the covering $X(2,12)\to
  X_1(12)$.  It has degree 2 and is branched over 4 cusps of the
  rational curve $X_1(12)$, so $X(2,12)$ has genus 1.
\end{proof}

\begin{rem}
  There are at least two other ways to bound the prime-to-$p$ torsion
  on $E$, both perhaps better for generalizations to higher genus.
  One uses an integrality result, namely
  Proposition~\ref{prop:integrality} below.  Another proceeds by
  considering the monodromy group of the finite-coefficient cohomology
  of the family $\EE\to\P^1$ associated to $E/\Fpbar(u)$.
\end{rem}

\section{The N\'eron model}\label{s:Neron}
Fix a positive integer $d$, let $F=\Fq(t^{1/d})=\Fq(u)$, and consider $E$ over
$F$.  In order to compute heights on $E(K_d)$ we will need a nice
model for $E$ over $\P^1$, specifically, a smooth proper surface $\EE_d$
over $\Fq$ equipped with a relatively minimal, generically smooth
morphism $\pi:\EE_d\to\P^1$ whose generic fiber is $E$.  

Using Tate's algorithm \cite{Tate75} we find that $E$ has reduction
type $I_{2d}$ at $u=0$ and type $I_2$ at the places dividing $u^d-1$.
At $u=\infty$, the reduction type is $I_{2d}$ if $d$ is even and
$I_{2d}^*$ if $d$ is odd.  At all other places, $E$ has good
reduction.  The reduction at $u=0$ and at $\infty$ when $d$ is even is
split multiplicative, while the reduction at the $d$-th roots of unity
is split if and only if $4$ divides $q-1$.

We only need the N\'eron model when $d=p^f+1$, which is even, so for
the rest of this section, \emph{we assume that $d$ is even\/}.  We
also assume for convenience of exposition that $d$ divides $q-1$, i.e.,
that $\Fq$ contains the $d$-th roots of unity.

We start by constructing a fibered projective surface $\WW_d\to\P^1$
whose fibers are plane cubics which are the reductions of $E$ at
places of $F$.  To that end, let $\WW_{d,0}$ be the closed subset of
$\P^3\times\A^1$ (with homogeneous coordinates $X$, $Y$, $Z$, on
$\P^3$ and coordinate $u$ on $\A^1$) defined by
$$Y^2Z-X(X+Z)(X+u^dZ)=0$$
and let $W^0_{d,0}$ be the open subset where $u\neq0$.
Also, let $\WW_{d,1}$ be the closed subset of
$\P^3\times\A^1$ (with homogeneous coordinates $X'$, $Y'$, $Z'$, on
$\P^3$ and coordinate $v$ on $\A^1$) defined by
$$Y^{\prime2}Z'-X'(X'+Z')(X'+v^dZ')=0$$
and let $W^0_{d,1}$ be the open subset where $v\neq0$.
We glue $W_{d,0}$ and $W_{d,1}$ along their open subsets $W^0_{d,0}$
and $W^0_{d,1}$ via the identification
$$([X',Y',Z'],v) = ([u^{-d}X,u^{-3d/2}Y,Z],u^{-1}).$$
We write $\WW_d$ for the resulting surface.  It is equipped with a
canonical morphism $\WW_d\to\P^1$ which is projective and generically
smooth and whose fibers are plane cubics.

Applying the Jacobian criterion, we find that $\WW_d$ is regular away
from $d+2$ points.  In the coordinates above, they are located at $u=0$,
$[X,Y,Z]=[0,0,1]$; $u\in\mu_d$, $[X,Y,Z]=[-1,0,1]$; and at $v=0$,
$[X',Y',Z']=[0,0,1]$.  The N\'eron model $\EE_d\to\P^1$ is obtained by
blowing up each of these singular points so that $\EE_d$ is smooth and
relatively minimal over $\P^1$. 

We first focus attention on the fiber over $u=0$.  The blow-ups
required to desingularize this fiber were discussed in some detail in
\cite[Lecture~3, \S3]{Ulmer11}, and we will not repeat all the details
here.  However, to establish notation to be used later, we will
quickly review what is needed.  We work in the affine open subset of
$\WW_{d,0}$ where $Z\neq0$ and set $x=X/Z$, $y=Y/Z$, so that our
surface is given by
$$y^2-x(x+1)(x+u^d)$$
in the affine 3-space with coordinate $x$, $y$, $u$.  The interesting
patch of the blow-up at $x=y=u=0$ has coordinates $x_1$, $y_1$, $u$
satisfying $x=ux_1$, $y=uy_1$, and in those coordinates, the proper
transform of our surface has equation
$$y_1^2-x_1(ux_1+1)(x_1+u^{d-1})$$
and the exceptional divisor is the reducible curve $u=y_1^2-x_1^2=0$.
The point $u=x_1=y_1=0$ is the unique singular point of the proper
transform, and we continue by blowing it up.  The new interesting
patch has coordinates $x_2$, $y_2$, $u$, with $x_1=ux_2$ and
$y_1=uy_2$; two new components are introduced; and there is a unique
singular point, which we again blow up.  At each of the first $d-1$
blow-ups, two new components are introduced.  At the $d$-th blow-up,
one component is introduced and we arrive at a smooth surface mapping
to $\P^1$ with fiber over $u=0$ consisting of a chain of $2d$ $\P^1$s
meeting transversally.  We label the components cyclically by
$i\in\Z/2d\Z$ so that the original fiber over $u=0$ is labeled 0 and
the component introduced in the first blow up with $u=y_1-x_1=0$ is
labeled 1.  (Thus the component $u=y_1+x_1=0$ is labeled $2d-1$.)

The situation over $v=0$ is essentially identical to that over $u=0$,
and we employ similar notation and conventions there.

Now consider the fiber over $u=\zeta^i_d$.  We set $x'=x+1=X/Z+1$,
$u'=u-\zeta^i_d$, and $f=(u^d-1)/(u-\zeta^i_d)$, so that $f$ is a unit
in the local ring at $u'=0$.  An affine patch of $\WW_d$ is
given by
$$y^2-(x'-1)x'(x'+u'f)=0.$$
The singular point is at $u'=x'=y=0$ and the tangent cone is the
irreducible curve with projective equation $y^2+x^{\prime2}+x'u'f=0$.
Blowing up once introduces a single component and the proper transform
of $\WW_d$ is regular in a neighborhood of $u'=0$.  This is consistent
with our earlier observation that the fibers over the points with
$u\in\mu_d$ are of type $I_2$.

The end result is a fibered surface $\pi:\EE_d\to\P^1$ such that
$\EE_d$ is regular and proper, and $\pi$ is generically smooth over
$\P^1$ with generic fiber $E/F$.  

The following result follows from our construction of $\EE_d$; see
\cite[Lecture 3, \S2]{Ulmer11} and especially Exercise~2.2 in that section.
The quantity $\delta$ defined here is sometimes called the
\emph{height} of $\EE_d$.

\begin{lemma}\label{lemma:height}
  Assume that $d$ is even.  Let $\pi:\EE_d\to\P^1$ be the model of
  $E/\Fq(t^{1/d})$ just constructed, with zero section $s:\P^1\to\EE_d$,
  and let $\omega$ be the invertible sheaf
  $s^*\Omega^1_{\EE_d/\P^1}\cong(R^1\pi_*\O_{\EE_d})^{-1}$.  Then $\delta$,
  the degree of $\omega$, is $d/2$.
\end{lemma}

\section{Heights}\label{s:heights}
Our next task is to compute the height pairing on the subgroup of
$E(F)$ generated by the points $P_i$.  This will give another proof
that the subgroup they generate has rank at least $d-2$ and it will give us
some control on the index in $E(F)$ of this subgroup. 

Throughout this section, we assume that $F=\Fq(t^{1/d})$ where
$d=p^f+1$ and $d$ divides $q-1$, and that $\pi:\EE_d\to\P^1$ is the
N\'eron model of $E/F$ discussed in the preceding section.

Both the height pairing and the discussion of discriminants in the
next section are based heavily on the Shioda-Tate isomorphism. We
refer to \cite[Chap.~4]{UlmerCRM} for a discussion of this in a
somewhat more general context.  In particular, the definition of the
height pairing on $E(F)$ and the proof that it is closely related to
the canonical height was discussed in \cite[Prop.~4.3.1]{UlmerCRM}.

We quickly review the definition of the height pairing here.  For a
point $P\in E(F)$, we write $[P]$ for the corresponding curve in
$\EE_d$, namely the image of the section $\P^1\to\EE_d$ whose generic
fiber is $P$.  For each $P$ there is a unique $\Q$-linear combination
$D_P$ of non-identity components of reducible fibers of $\pi$ such
that the divisor $[P]-[0]+D_P$ is orthogonal to all components of
fibers of $\pi$ under the intersection pairing on $\EE_d$.  To compute
$D_P$, one need only determine the component of the fiber through
which $[P]$ passes at each reducible fiber.  The height pairing is
then defined by
\begin{equation}\label{eq:height}
  \langle P,Q\rangle=-\left([P]-[0]+D_P\right)\cdot\left([Q]-[0]\right).
\end{equation}
Here the dot signifies the intersection pairing on $\EE_d$.  

This pairing is known to be symmetric and positive definite modulo
torsion. It makes $E(F)$ into a lattice in the real vector space
$E(F)\tensor_\Z\R$.  The N\'eron-Tate canonical height is $\log q$
times this height.  

The ``correction factor'' $D_P.[Q]$ is symmetric
in $P$ and $Q$ and depends only on the components through which $[P]$
and $[Q]$ pass at each reducible fiber.  A table of values is given in
\cite[1.19]{CoxZucker79}.

For any elliptic surface $\pi:\EE\to\P^1$, the self intersection of any section
is equal to $-\delta$ where $\delta$ is the degree of the invertible
sheaf $\omega=(R^1\pi_*\O_\EE)^{-1}$.  We computed in
Lemma~\ref{lemma:height} that for $\EE_d$ we have $\delta=d/2$.

We apply the above to the elliptic surface $\EE_d\to\P^1$ associated
to $E$ over $F$ and the points $P_i$ and obtain the following.

\begin{thm}\label{thm:heights}
  The height pairing~\eqref{eq:height} of the points $P_i$
  \textup{(}$i\in\Z/d\Z$\textup{)} on $E$ given by~\eqref{eq:points} is
\begin{equation*}
  \langle P_i,P_j\rangle =
\begin{cases}
\frac{(d-1)(d-2)}{2d}&\text{if $i=j$}\\
\\
\frac{(1-d)}{d}&\text{if $i-j$ is even and $\neq0$}\\
\\
0&\text{if $i-j$ is odd}
\end{cases}
\end{equation*}
\end{thm}

\begin{proof}
  Throughout, we conflate the point $P_i$ with the section of
  $\EE_d\to\P^1$ determined by $P_i$.

  Let us first examine the reduction of
  $P_i=(\zeta_d^iu,\zeta_d^iu(\zeta_d^iu+1)^{d/2})$ in the various
  reducible fibers.  At $u=0$ it is clear that $P_i$ reduces to the
  singular point in the fiber of $\WW_d$.  In the first blow up, $P_i$
  has coordinates
  $(x_1,y_1)=(\zeta_d^i,\zeta_d^i(\zeta_d^iu+1)^{d/2})$ and so reduces
  to a smooth point on the component labeled 1.  At $u=\zeta_d^j$, we
  see that $P_i$ reduces to a point on the identity component if and
  only if $\zeta_d^{i+j}\neq-1$, i.e., if and only if $i+j\not\equiv
  d/2\pmod d$.

  To examine the reduction at $u=\infty$ (i.e., over $v=0$ in the
  coordinates used to construct $\WW_d$ and $\EE_d$), we introduce
  $x'=X'/Z'$ and $y'=Y'/Z'$.  In these coordinates, we have that
  $x'(P_i)=\zeta_d^iv^{d-1}$ and
  $y'(P_i)=\zeta_d^iv^{d-1}(\zeta_d^i+v)^{d/2}$.  Thus in the model
  $\WW_d$, $P_i$ reduces to the singular point in the fiber over $v=0$
  and we have to follow its behavior in the blow ups used to
  construct $\EE_d$.  Using notation analogous to that in the last
  section, we see that for $1\le a\le d-1$,
  $x'_a(P_i)=\zeta_d^iv^{d-1-a}$ and so $P_i$ continues to reduce to a
  singular point until the $(d-1)$st blow up.  At that stage, $P_i$ has
  coordinates
$$(x'_{d-1},y'_{d-1})=(\zeta_d^i,\zeta_d^i(\zeta_d^i+v)^{d/2}).$$
If $i$ is even, then $(\zeta_d^i)^{d/2}=1$ and so $P_i$ reduces to the
component labeled $d-1$, whereas if $i$ is odd, then
$(\zeta_d^i)^{d/2}=-1$ and so $P_i$ reduces to the component labeled
$d+1$.

This calculation allows us to determine the ``correction factors''
$D_{P_i}.[P_j]$, using the table \cite[1.19]{CoxZucker79}.  At $u=0$,
both points reduce to the component labeled 1, so the contribution is
$(2d-1)/(2d)$.  At places with $u\in\mu_d$, if $i\neq j$ then at least
one of $P_i$ or $P_j$ lands on the identity component, so the
contribution is always 0.  If $i=j$, then at exactly one root of unity
(namely at $u=\zeta_d^{d/2-i}=-\zeta_d^{-i}$), $P_i$ lands on the
non-identity component, and we get a local contribution of $1/2$.
Finally, at $u=\infty$, $P_i$ and $P_j$ land on components $d\pm1$ and
we find a local contribution of 
$$\begin{cases}
\frac{(d-1)(d+1)}{2d}&\text{if $i-j$ is even}\\
\frac{(d-1)^2}{2d}&\text{if $i-j$ is odd.}
\end{cases}
$$

The above calculations also show that the intersection
numbers $[P_i].[P_j]$ ($i\neq j$) and $[P_i].[O]$ are 0.  Indeed, it
is evident from inspecting the $x$-coordinates that $P_i$ and $P_j$ do
not meet except possibly at $u=0$ or $u=\infty$.  But at $u=0$, $P_i$
reduces to $(x_1,y_1)=(\zeta_d^i,\zeta_d^i)$;
thus $P_i$ and $P_j$ do not meet in that fiber. At $u=\infty$, $P_i$ reduces to
$(\zeta_d^i,\pm\zeta_d^i)$, and so again $P_i$ and $P_j$ do not meet in
that fiber.  To see that $[P_i].[O]=0$ we note that the coordinates of
$P_i$ are polynomials in $u$, so $P_i$ cannot meet the zero
section except at $u=\infty$.  But near $u=\infty$ (i.e., near $v=0$) we
just saw that the coordinates of $P_i$ are polynomials in $v$ and so
$P_i$ does not meet $O$ in that fiber either.
Thus we have that 
$$(P_i-O).(P_j-O)=\begin{cases}
P_i^2+O^2=-d&\text{if $i=j$}\\
O^2=-d/2&\text{if $i\neq j$}
\end{cases}
$$

Combining the ``geometric'' contribution $-(P_i-O).(P_j-O)$ with the
``correction factor'' $-D_{P_i}.P_j$ leads, after a small bout of
arithmetic, to the stated result.  This completes the proof of the
theorem.
\end{proof}

We could have simplified the calculation slightly by noting that since
$E$ is defined over $K_1=\Fp(t)$ and $\gal(F/K_1)$ permutes the $P_i$
cyclically, the pairing $\langle P_i,P_j\rangle$ depends only on
$i-j$ modulo $d$.

\begin{rem}
We write $V_d$ for the subgroup of $E(K_d)$ generated by the $P_i$.
It is easy to see that the lattice $V_d$ modulo torsion is isomorphic
to the direct sum of two copies of a scaling of the lattice
$A^*_{d/2-1}$ where $A^*_n$ is the dual of the $A_n$ lattice, as
in~\cite[4.6.6]{ConwaySloaneSPLG}.  The points $P_i$ are root
vectors; perhaps this explains why it is so easy to write them down.
\end{rem}

We note that the height calculation provides a second proof that $V_d$
has rank $d-2$.

We finish this section by tying up a loose end from
Section~\ref{s:points}.

\begin{prop}\label{prop:trace-to-d=2}
  Suppose that $d=p^f+1$.  Then for $j\in\Z/2\Z$, we have the
  following equality in $E(K_d)$:
$$\sum_{\substack{i\in\Z/d\Z\\i\equiv j+1+d/2\pmod2}}P^{(d)}_i=P^{(2)}_j.$$
\end{prop}

\begin{proof}
  We saw in Section~\ref{s:points} that the sum on the left hand side
  is a 4-torsion point.  Applying the map $\nu$ of
  Proposition~\ref{prop:Sel}, we see that both sides of the equality
  above map to either $(1,0,\dots,1,0)$ or $(0,1,\dots,0,1)$ in
  $(\Z/2\Z)^d$ (depending on the parity of $j$).  This shows that the
  equality above holds in $E(K_d)/2E(K_d)$ and therefore holds in
  $E(K_d)$ up to a sign.  (The kernel of $E(K_d)_{tor}\to
  E(K_d)/2E(K_d)$ is generated by $Q_0$, and we have
  $P^{(2)}_j+Q_0=-P^{(2)}_j$.)  To pin down the sign, we consider the
  images of the two sides in the component group of $E$ at $u=0$.  We
  saw in the proof of Theorem~\ref{thm:heights} that $P^{(d)}_i$
  reduces to the component labeled 1, so the sum reduces to the
  component labeled $d/2$.  To finish, one checks that $P^{(2)}_j$
  reduces to component $d/2$ while $-P^{(2)}_j$ reduces to component
  $2d-d/2$.  This completes the proof of the Proposition.
\end{proof}

\section{The index $[E(F):V_d]$}\label{s:index}
We start by proving a result that is probably well-known to experts.
There is a closely related discussion in a more general context in
\cite[\S5]{Gordon79}), but I do not know a convenient reference for
the result we need.

Suppose that $F=\Fq(\CC)$ is the function field of a curve over $\Fq$
and $E/F$ is an elliptic curve with N\'eron model $\pi:\EE\to \CC$.
For each place $v$ if $F$, let $f_v$ be the number of irreducible
components in the fiber of $\pi$ over $v$.  One of these components
meets the zero section and we let $N_v$ be the subgroup of $\NS(\EE)$
generated by the other components of the fiber.  This is known to be a
free abelian group of rank $f_v-1$ (e.g., \cite[Lecture 2,
8.6]{Ulmer11}).  We write $\disc_v$ for the discriminant of the
intersection form on $\EE$ restricted to $N_v$. (If $f_v=1$, we set
$\disc_v=1$.)  We write $\langle\cdot,\cdot\rangle$ for the height
pairing on $E(F)$ discussed in the previous section (so that
$\langle\cdot,\cdot\rangle \log q$ is the N\'eron-Tate canonical
height).  Let $R_1,\dots,R_r$ be a basis for $E(F)/tor$ and let
$$R(E/F)=\left|\det(\langle R_i,R_j\rangle)_{i,j=1,\dots,r}\right|.$$
This is the usual regulator of $E/F$ divided by $(\log q)^r$.

The result we have in mind gives an \emph{a priori} bound on the denominator
of the rational number $R(E/F)$.

\begin{prop}\label{prop:integrality}
The rational number 
$$\frac{R(E/F)\prod_v\disc_v}{|E(F)_{tor}|^2}$$
is an integer.
\end{prop}

\begin{proof}
  Let $L^1\NS(\EE)$ be the subgroup of the N\'eron-Severi group
  $\NS(\EE)$ generated by classes of divisors which meet the generic
  fiber in a divisors of degree 0, and let $L^2\NS(\EE)$ be the
  subgroup generated by components of fibers of $\pi$.  The Shioda-Tate
  isomorphism says that restriction to the generic fiber induces an
  isomorphism $L^1\NS(\EE)/L^2\NS(\EE)\to E(F)$.  See
  \cite[Lecture 3, \S5]{Ulmer11} for a discussion in our context, as well as
  \cite[Chap.~4]{UlmerCRM} for a proof of a more general result and
  several results used below.

  The intersection form restricted to $L^1\NS(\EE)$ is degenerate---by
  definition it is orthogonal to $F$, the class of a fiber of $\pi$.
  We define $\overline L^1$ as $L^1\NS(\EE)$ modulo $\Z F$.  Since
  $\NS(\EE)$ is torsion-free \cite[VII.1.2]{MirandaBTES}
  and $F$ is primitive, $\overline L^1$ is torsion free.  The
  intersection form induces a non-degenerate integral pairing on
  $\overline L^1$.  For any subgroup $L\subset \overline L^1$, we
  write $\disc(L)$ for the discriminant of the intersection pairing
  restricted to $L$.  (By discriminant, we mean the \emph{absolute
    value} of the matrix of pairings.)

We write $\overline L^2$ for the image of $L^2\NS(\EE)$ in $\overline
L^1$.  We have an orthogonal direct sum decomposition
$$\overline L^2=\prod_vN_v.$$

  Recall that $R_1,\dots,R_r$ form a basis for $E(F)/tor$.  For each
  $i$ we define the class
$$\widetilde R_i=[R_i]-[O]$$
(the class of the section of $\pi$ attached to $R_i$ minus the zero
section).  By the Shioda-Tate isomorphism, the subgroup $N$ of
$\overline L^1$ generated by $\widetilde R_1,\dots,\widetilde R_r$ and
$\overline L^2$ has index equal to $|E(F)_{tor}|$. Thus we have
$$\disc(N)=|E(F)_{tor}|^2\disc(\overline L^1).$$

Let $t$ be the least common multiple of the exponents of the component
groups over all place of $F$, so that for each $i$, the section
associated to $tR_i$ passes through the identity component at each
place.  Let $\widetilde{tR_i}$ be the class of $[tR_i]-[O]$, and let $N'$
be the subgroup of $\overline L^1$ generated by $\widetilde{tR_1},\dots,
\widetilde{tR}_r$ and $\overline L^2$.  The index of $N'$ in $N$ is clearly $t^r$ and
so $\disc(N')=t^{2r}\disc(N)$.

On the other hand, since each $tR_i$ passes through the identity
component at each place $v$, the $\widetilde{tR_i}$ are orthogonal to
$\overline L^2$.  Thus, we have an orthogonal direct sum decomposition
$$N'=\left\{\Z\widetilde{tR_1}+\cdots+\Z\widetilde{tR_r}\right\}
\bigoplus\overline L^2.$$
Also, in the height pairing $\langle tR_i,tR_j\rangle$, no
``correction factors'' are required, i.e., we have
$$\langle tR_i,tR_j\rangle=-\widetilde{tR_i}.\widetilde{tR_j}.$$

Taking discriminants, we find
\begin{align*}
\disc(N')&=t^{2r}R(E/F)\disc(\overline L^2)\\
&=t^{2r}R(E/F)\prod_v\disc N_v.
\end{align*}

But we saw above that
$$\disc(N')=t^{2r}\disc(N)=t^{2r}|E(F)_{tor}|^2\disc(\overline L^1).$$
Thus
$$\frac{R(E/F)\prod_v\disc_v}{|E(F)_{tor}|^2}=\disc(\overline L^1)$$ 
which is an integer.
\end{proof}

We now apply the proposition in the context of the Legendre curve $E$.
Let $F=\Fq(t^{1/d})$ with $d=p^f+1$ for some positive integer $f$ and
assume that $d$ divides $q-1$.  Let $V_d$ be the subgroup of
$E(F)$ generated by the points $P_0,\dots,P_{d-1}$ of
Section~\ref{s:points}.  We saw in Section~\ref{s:upper-bound} that
$V_d$ has finite index in $E(F)$.  Using the proposition, we find the
following \emph{a priori} bound on the index.

\begin{thm}\label{thm:index}
The index $[E(F):V_d]$ divides $p^{f(d-2)/2}=p^{f(p^f-1)/2}$.  In
particular, it is a power of $p$.
\end{thm}

\begin{proof}
  The lattice $A_n^*$ has discriminant $(n+1)^{n-1}$
  \cite[4.6.6]{ConwaySloaneSPLG} and so Theorem~\ref{thm:heights}
  shows that the discriminant of $V_d$ modulo torsion is
  $2^{4-d}(d-1)^{d-2}d^{-2}$.  Thus, we have that the regulator is
\begin{equation}\label{eq:reg}
R(E/F)=\frac{2^{4-d}(d-1)^{d-2}d^{-2}}{[E(K_d):V_d]^{2}}.
\end{equation}

The local factors
$\disc_v$ at $u=0$ and $u=\infty$ are $2d$, at the $d$ places
where $u^d=1$ they are $2$, and elsewhere they are 1.  (These are
also the orders of the component groups at these places.  This is of
course not an accident---see \cite[Thm~5.2ff]{Gordon79}.)
By Proposition~\ref{prop:torsion}, $|E(F)_{tor}|=8$.
Unwinding the integrality statement in the proposition, it follows
that
$$\frac{p^{f(d-2)}}{[E(F):V_d]^2}$$
is an integer.  This completes the proof of the theorem.
\end{proof}

\begin{rems}
Let $f$ and $f'$ be positive integers such that $f'$ divides $f$ and
$f/f'$ is odd.  Let $d=p^f+1$ and $d'=p^{f'}+1$ and note that $d'$
divides $d$ so that $K_{d'}\subset K_d$.  
\begin{enumerate}
\item Writing $\langle\cdot,\cdot\rangle_{d}$ and
  $\langle\cdot,\cdot\rangle_{d'}$ for the height pairings on $E(K_d)$
  and $E(K_{d'})$ respectively, we have that 
$$\langle P^{(d')}_0,P^{(d')}_0\rangle_{d'}=\frac{(d'-1)(d'-2)}{2d'}$$
and so
$$\langle P^{(d')}_0,P^{(d')}_0\rangle_{d}=\left(\frac{d}{d'}\right)\frac{(d'-1)(d'-2)}{2d'}.$$
On the other hand, the root vectors $P^{(d)}_i$ of $V_d$ satisfy
$$\langle P^{(d)}_i,P^{(d)}_i\rangle_{d}=\frac{(d-1)(d-2)}{2d}.$$
The ratio 
$$\frac{\langle P^{(d)}_i,P^{(d)}_i\rangle_{d}}{\langle P^{(d')}_0,P^{(d')}_0\rangle_{d}}=
\frac{(d-1)}{(d-(d/d'))}\frac{(d-2)}{(d-(2d/d'))}$$
is $>1$ and so $P^{(d')}_0$ is not contained in $V_d$.  This proves
that if $f$ is divisible by an odd prime, then $V_d\neq E(K_d)$.
\item In a future publication, we will compute the quotient group $E(K_d)/V_d$ as a
  module over $\Z_p[\gal(K_d/K_1)]$ and show that it is non-trivial
  precisely when $f>2$. In general it is not generated by points in
  $V_{d'}$ as $d'$ runs through divisors of $d$ of the form $d'=p^{f'}+1$.
\item With more work, one may compute the height pairings
  $\langle P^{(d')}_j,P^{(d)}_i\rangle_d$ and show that for every
$j\in\Z/d'\Z$, we have the following equality in $E(K_d)$ modulo
torsion:
$$\sum_{\substack{i\in\Z/d\Z\\i\equiv j\pmod{d'}}}P^{(d)}_i=-p^{(f-f')/2}P^{(d')}_{j.}$$
(Note the surprising sign.) 
\end{enumerate}
\end{rems}

\section{The Tate-Shafarevich group}\label{s:sha}
In this section we use the Birch and Swinnerton-Dyer formula to relate
the order of the Tate-Shafarevich group to the index $[E(F):V_d]$ of
the preceding section.  We refer to \cite{Ulmer11} and \cite{UlmerCRM}
for background on $L$-functions and the BSD conjecture.

First we calculate the $L$-function of $E$ over $\Fq(t^{1/d})$ when
$d=p^f+1$.

\begin{prop}
Let $F=\Fq(t)(t^{1/d})$ where $d=p^f+1$ and $d$ divides $q-1$.  Then
the Hasse-Weil $L$-function of $E/F$ is
$$L(E/F,s)=(1-q^{1-s})^{d-2}.$$
Moreover, the Birch and Swinnerton-Dyer conjecture holds for $E/F$.
\end{prop}

\begin{proof}
  By \cite[Lecture 1, Theorem~9.3]{Ulmer11}, the $L$-function $L(E/F,s)$ is a
  polynomial in $q^{-s}$ of the form $\prod_i(1-\alpha_iq^{-s})$ where
  the $\alpha_i$ are Weil integers of size $q$.  Its degree is
  $\deg\n-4$ where $\n$ is the conductor of $E$.

  We saw in Section~\ref{s:Neron} that the curve $E$ has
  multiplicative reduction at $d+2$ places of $F$ and good reduction
  elsewhere.  It follows that its conductor has degree $d+2$.  Thus
  the degree of the $L$-function as a polynomial in $q^{-s}$ is
  $d-2$.  Note that the order of vanishing of $L(E/F,s)$ at $s=1$ is
  precisely the number of times $q$ appears as an $\alpha_i$.

  Now we recall (e.g., from \cite[Theorem~6.3.1]{UlmerCRM}) that we
  have an inequality 
$$\ord_{s=1}L(E/F,s)\ge\rk E(F).$$  
Since the degree of $L(E/F,s)$ is $d-2$, the order of vanishing is
certainly at most $d-2$.  On the other hand, we know from
Corollary~\ref{cor:rks} that under the hypotheses of the proposition,
the rank of $E(F)$ is $d-2$.  Thus all of the $\alpha_i$ must be equal
to $q$, i.e., the $L$-function is as we have claimed. Moreover, we
have the equality
$$\ord_{s=1}L(E/F,s)=\rk E(F)$$
i.e., the basic BSD conjecture holds for $E/F$.  
\end{proof}

Using the refined BSD conjecture (on the leading term of the
$L$-function) we get the following connection between the
Tate-Shafarevich group $\sha(E/F)$ and the index.

\begin{cor}\label{cor:sha}
Suppose that $d=p^f+1$ and let $K_d=\Fp(\mu_d,t^{1/d})$.  Recall that
$V_d$ is the subgroup of $E(K_d)$ generated by the points
$P_0,\dots,P_{d-1}$.  We have that
$\sha(E/K_d)$ is finite and
$$|\sha(E/K_d)|=[E(K_d):V_d]^2.$$
More generally, if $F=\Fq(t^{1/d})$ with $d$ dividing $q-1$, we have
$$|\sha(E/F)|=[E(F):V_d]^2(q/p^{2f})^{(p^f-1)/2}.$$
In particular, $\sha(E/F)$ is a $p$-group.
\end{cor}

This result can be viewed as a class number formula, analogous to that
relating the class number to the index of cyclotomic units or the
index of a Heegner point to the order of a Tate-Shafarevich group
(cf. \cite[Thm.~8.2]{WashingtonCF} and
\cite[Conj.~V.2.2]{GrossZagier86}).

\begin{proof}
It is known (see,
e.g., \cite[Theorem~6.3.1]{UlmerCRM}) that the basic BSD conjecture
implies the refined BSD conjecture.  For $E/F$, writing $r=d-2$ for
the rank, we have that the leading coefficient is 
$$(1/r!)L^{(r)}(E/F,s)=(\log q)^r.$$  
On the other hand, the regulator $\reg(E/F)$ (computed using the
N\'eron-Tate canonical height) is equal to $(\log q)^rR(E/F)$ where
$R(E/F)$ is the determinant appearing in Section~\ref{s:heights}.
Thus the refined BSD conjecture leads to the equality
$$1=\frac{|\sha(E/F)|R(E/F)\tau(E/F)}{|E(F)_{tor}|^2}$$
where $\tau(E/F)$ is the Tamagawa number.  

To unwind this, recall (e.g., from \cite[\S6.2.3]{UlmerCRM}) that 
$$\tau(E/F)=q^{1-\delta}\prod_vc_v$$
where $\delta$ is the integer appearing in Lemma~\ref{lemma:height}
and $c_v$ is the order of the component group of the N\'eron model at
$v$.  Explicitly, we find that
$$\tau(E/F)=q^{1-d/2}2^d(2d)^2=q^{1-d/2}2^{d+2}d^2.$$

Using the value of $R(E/F)$ from Equation~\eqref{eq:reg}
and the order of the torsion subgroup from
Proposition~\ref{prop:torsion}, a short calculation leads to
$$|\sha(E/F)|=[E(F):V_d]^2(q/p^{2f})^{(p^f-1)/2}.$$
Since we proved in Theorem~\ref{thm:index} that the index $[E(F):V_d]$
is a power of $p$, so is the order of $\sha(E/F)$.  This completes the
proof of the corollary.
\end{proof}

\begin{rem}
Note that the upper bound on $[E(F):V_d]$ of Theorem~\ref{thm:index}
gives an upper bound on the order of $\sha(E/F)$: it has order at most
$q^{(d-2)/2}$.  By the formula of the corollary, its order is bounded
below by $(q/p^{2f})^{(d-2)/2}$.  This lower bound suggests that there
might be a construction of homogenous spaces for $E$ with parameters
with values in $\Fq$ such that the homogenous space is trivial when the
parameters lie in $\Fp(\mu_d)$.
\end{rem}

\begin{rem}
  We computed the $L$-function when $d=p^f+1$ by using the BSD
  inequality and a lower bound on the rank.  We could also deduce the
  $L$-function from \cite[4.7]{Ulmer07b} or from the Artin formalism
  and a root number calculation without knowing anything about points
  on $E(F)$.
\end{rem}

\section{Connection with Berger's construction}\label{s:Berger}
In this section, we let $k$ be any field of characteristic $\neq2$, and
we relate the Legendre curve over $k(t)$ to a curve studied in
\cite{Ulmer13a}.  This leads to strong consequences for the
arithmetic of $E$ at every layer of the tower of fields $k(t^{1/d})$.

Let $E'$ over $k(t')$ be the elliptic curve studied in Section~7 of
\cite{Ulmer13a}:
$$E':\quad y^{\prime2}+x'y'+t'y'=x^{\prime3}+t'x^{\prime2}.$$
We identify the fields $k(t)$ and $k(t')$ by setting $t'=t/16$ and
view $E'$ as an elliptic curve over $k(t)$.

\begin{lemma}
The elliptic curve $E'$ over $k(t)$ is 2-isogenous to the Legendre
curve~\eqref{eq:E}.
\end{lemma}

\begin{proof}
This follows from a simple computation.  On $E'$, changing coordinates
$$(x',y')\mapsto\left(\frac{x'-t}{16},\frac{y'-2x'}{64}\right)$$
leads to the equation
$$y^{\prime2}=x^{\prime3}+(4-2t)x^{\prime2}+t^2x'.$$
Dividing by the 2-torsion point $(0,0)$ yields the curve
$$y^2=x^3+(4t-8)x^2-16(t-1)x=x(x-4)(x+4(t-1)).$$
(See, e.g., \cite[Chap.~4, \S5]{HusemollerEC}.) 
Finally, the change of coordinates $(x,y)\mapsto(4x+4,8y)$ yield the
Legendre curve
$$y^2=x(x+1)(x+t).$$
\end{proof}

\begin{rem}  Let $\psi:E\to E'$ be the isogeny dual to the one of the
  lemma (which we have only specified up to a sign).  Writing
  $\psi^*$ for the pull back of functions,  we have
$$\psi^*(x')=\frac{x^2-t}{16(x+1)}$$
and
$$\psi^*(y')=\pm\left(\frac{y(x^2+2x+t)}{64(x+1)^2}-\frac{x(x+t)}{32(x+1)}\right).$$ 
  Under this isogeny, the point $P_0$ of this paper maps to the sum of
  $(-t',0)$ and the point $-P(u')$ in
  \cite[Theorem~8.1(3)]{Ulmer13a}.  In particular, the group $V_d$ of
  this paper is mapped to the group $V_d$ of
  \cite[Remark~8.10(3)]{Ulmer13a}.
\end{rem}

The lemma allows us to transfer strong results about $E'$ to $E$.

\begin{cor}\label{cor:BSD}
If $k$ is a field of characteristic zero, then the rank of
$E(k(t^{1/d}))$ is zero for all $d$.  If $k$ is a finite field of
characteristic $p>2$, then the Birch and Swinnerton-Dyer conjecture
holds for $E$ over $k(t^{1/d})$ for all $d$ prime to $p$.
\end{cor}

\begin{proof}
  The conclusions (rank zero or BSD) both ``descend'' under finite
  extensions of $k$, i.e., if they hold for $E$ over $k'(t^{1/d})$
  with $k'$ a finite extension of $k$, then they hold over
  $k(t^{1/d})$.  So for a fixed $d$ we may assume that $k$ contains a
  $d$-th root of $16$.  In this case, we may take $k(t^{1/d})$ and
  $k(t^{\prime 1/d})=k((t/16)^{1/d})$ to be the same field inside an
  algebraic closure of $k(t)$.  The lemma shows that $E$ and $E'$ are
  2-isogenous over this field and we proved the claims of the lemma
  for $E'$ in \cite[Section~7]{Ulmer13a}.  This shows that these
  claims also hold for $E$.
\end{proof}

\section{Summary and future directions}\label{s:summary}
For the convenience of the reader, we summarize the main results of
the paper.

\begin{thm}\label{thm:main}
Let $k$ be a field of characteristic $\neq2$ and let $k(t)$ be the
rational function field over $k$.  Let $E$ be the elliptic curve over
$k(t)$ given by
$$y^2=x(x+1)(x+t).$$
\begin{enumerate}
\item If the characteristic of $k$ is zero, then for all positive
  integers $d$, the rank of $E(k(t^{1/d}))$ is zero.
\item If $k$ is finite of characteristic $p$, then for all positive
  integers $d$ not divisible by $p$, the Birch and Swinnerton-Dyer
  conjecture holds for $E$ over $k(t^{1/d})$.  In particular, the
  Tate-Shafarevich group $\sha(E/k(t^{1/d}))$ is finite.
\item Suppose  that $k$ is finite of characteristic $p$ and
  cardinality $q$, and that $d$ divides $p^f+1$ for some $f$. Then
$$\rk E(k(t^{1/d}))=
\sum_{\substack{e|d\\e>2}}\frac{\varphi(e)}{o_q(e)}$$
where the sum is over divisors $e$ of $d$ greater than 2, $\varphi$ is
Euler's function, and $o_q(e)$ is the order of $q$ in
$(\Z/e\Z)^\times$.
\item Suppose further that $d=p^f+1$, $k=\Fp(\mu_d)$, and
  $K_d=k(t^{1/d})$.  Let $P_0,\dots,P_{d-1}$ be the points in $E(K_d)$
  defined by Equation~\eqref{eq:points}, and let $V_d$ be the subgroup
  of $E(K_d)$ generated by the $P_i$.  Then $V_d$ has rank $d-2$, its
  index in $E(K_d)$ is a power of $p$, and we have
$$[E(K_d):V_d]^2=|\sha(E/K_d)|.$$
\end{enumerate}
\end{thm}

Parts (1) and (2) of the Theorem were proven in
Corollary~\ref{cor:BSD}, part (3) is Corollary~\ref{cor:rks} (4), and
part (4) is the conjunction of Corollary~\ref{cor:rks} (4),
Theorem~\ref{thm:index}, and Corollary~\ref{cor:sha}.

\begin{rem}
  The claim in part (2) of the theorem in fact holds for all $d$.
  Indeed, if $d=pd'$, then $E/k(t^{1/d'})$ is isogenous to
  $E^{(p)}/k(t^{1/d'})$ and the latter is isomorphic to $E/k(t^{1/d})$
  after identifying $k(t^{1/d})\cong k(u)$ with $k(u')\cong
  k(t^{1/d'})$.  Since the truth of the BSD conjecture is invariant
  under isogeny, we can reduce to the case where $d$ is not divisible
  by $p$.
\end{rem}

For $k$ finite and $d$ dividing $p^f+1$ for some $f$, we have
calculated the rank of $E(k(t^{1/d}))$ now in three different ways:
via cohomology in \cite{Ulmer07b}, using Berger's construction and
endomorphisms of Jacobians in \cite{Ulmer13a}, and ``by hand'' in
this paper.  

The Theorem suggests (at least) two further directions for research
that we will address in future publications.  First, what is the rank
of $E/K_d$ for values of $d$ that do not divide $p^f+1$ for any $f$?
Surprisingly, there are many other ``interesting'' values of $d$
\cite{CHU} and, in a sense made precise in \cite{PomeranceUlmer},
those dividing $p^f+1$ are less numerous than the others.

Second, what can be said about the index $[E(K_d):V_d]$ and the
Tate-Shafarevich group, say when $d=p^f+1$?  In a future publication
we will use flat cohomology and crystalline cohomology to show that
the index is 1 if and only if $f\le2$ and to compute the quotient
$E(K_d)/V_d$ and the Tate-Shafarevich group $\sha(E/K_d)$ as modules
for the group ring $\Zp[\gal(K_d/K_1)]$.

\bibliography{database}{}
\bibliographystyle{alpha}

\end{document}

%% file: formatting.tex
\numberwithin{subsection}{section}

\numberwithin{equation}{section}

\theoremstyle{plain}
\newtheorem{thm}[equation]{Theorem}
\newtheorem{prop}[equation]{Proposition}
\newtheorem{cor}[equation]{Corollary}
\newtheorem{lemma}[equation]{Lemma}

\theoremstyle{definition}

\theoremstyle{remark}
\newtheorem{rem}[equation]{Remark}
\newtheorem{rems}[equation]{Remarks}

\newtheorem{rem-exer}[equation]{Remark/Exercise}
\newtheorem{rem-exers}[equation]{Remark/Exercises}

\newtheorem{exs}[equation]{Examples}

%% file: macros.tex
\usepackage[OT2,T1]{fontenc}
\DeclareSymbolFont{cyrletters}{OT2}{wncyr}{m}{n}
\DeclareMathSymbol{\sha}{\mathalpha}{cyrletters}{"58}

\newcommand{\CC}{\mathcal{C}}

\newcommand{\WW}{{\mathcal{W}}}

\renewcommand{\O}{\mathcal{O}}

\newcommand{\EE}{\mathcal{E}}

\newcommand{\F}{\mathbb{F}}
\newcommand{\Fp}{{\mathbb{F}_p}}

\newcommand{\Fq}{{\mathbb{F}_q}}

\newcommand{\Fpbar}{{\overline{\mathbb{F}}_p}}

\newcommand{\Zp}{{\mathbb{Z}_p}}

\newcommand{\Z}{\mathbb{Z}}

\newcommand{\Q}{\mathbb{Q}}
\newcommand{\R}{\mathbb{R}}

\newcommand{\A}{\mathbb{A}}
\renewcommand{\P}{\mathbb{P}}

\newcommand{\n}{\mathfrak{n}}

\newcommand{\into}{\hookrightarrow}

\newcommand{\tensor}{\otimes}

\newcommand{\sdp}{{\rtimes}}

\newcommand{\labeledto}[1]{\overset{#1}{\to}}

\DeclareMathOperator{\im}{Im}

\DeclareMathOperator{\ord}{ord}
\DeclareMathOperator{\rk}{Rank}
\DeclareMathOperator{\reg}{Reg}
\DeclareMathOperator{\dvsr}{div}

\DeclareMathOperator{\NS}{NS}

\DeclareMathOperator{\gal}{Gal}

\DeclareMathOperator{\disc}{Disc}

\DeclareMathOperator{\Sel}{Sel}
\DeclareMathOperator{\Fr}{Fr}

\def\clap#1{\hbox to 0pt{\hss#1\hss}}